\documentclass[a4paper,11pt, reqno]{amsart}
\usepackage[all]{xy}
\usepackage{amsmath, amsthm, amsfonts, amssymb} 
\newtheorem{theo}{Theorem}[section]

\newtheorem{df}[theo]{Definition}

\newtheorem*{ack}{Acknowledgment}

\newtheorem{lem}[theo]{Lemma}

\newtheorem{nota}[theo]{Notation}

\newtheorem{prop}[theo]{Proposition}

\newtheorem{rem}[theo]{Remark}

\newtheorem{exam}[theo]{Example}

\newcommand{\R}{\mathbf{R}}
\newcommand{\C}{\mathbf{C}}
\newcommand{\Z}{\mathbf{Z}}
\newcommand{\F}{\mathbf{F}}

\newcommand{\N}{\mathbf{N}}
\newcommand{\Sd}{\operatorname{Sd}}
\newcommand{\pol}{\operatorname{pol}}

\newcommand{\Ad}{\operatorname{Ad}}
\newcommand{\id}{\operatorname{id}}
\newcommand{\Tr}{\operatorname{Tr}}

\begin{document}

\title[On some free products of von Neumann algebras]{On some free products of von Neumann algebras which are free Araki-Woods factors}

\author{Cyril Houdayer}

\address{Department of Mathematics \\
                 UCLA \\
                 Los Angeles \\
                 California \\
                 CA 90095}

\email{cyril@math.ucla.edu}




\begin{abstract}
We prove that certain free products of factors of type ${\rm I}$ and other von Neumann algebras with respect to nontracial, almost periodic states are almost periodic free Araki-Woods factors. In particular, they have the free absorption property and Connes' $\Sd$ invariant completely classifies these free products. For example, for $\lambda, \mu \in ]0, 1[$, we show that
$$(M_2(\C), \omega_{\lambda}) \ast 
(M_2(\C), \omega_{\mu})$$
is isomorphic to the free Araki-Woods factor whose $\Sd$ invariant is the subgroup of $\R^*_+$ generated by $\lambda$ and $\mu$. Our proofs are based on algebraic techniques and amalgamated free products. These results give some answers to questions of Dykema and Shlyakhtenko.
\end{abstract}

\maketitle

\section{Introduction}

In \cite{dykema96} and \cite{dykema94b}, Dykema investigated free products of finite dimensional and other von Neumann algebras with respect to nontracial faithful states. We are interested in free products of factors of type ${\rm I}$. In this respect, we recall Theorem $1$ of \cite{dykema96} in the particular case of factors of type ${\rm I}$ (see Proposition $7.3$ of \cite{dykema96} for a precise statement in full generality).
\begin{theo}[Dykema, \cite{dykema96}]\label{dykema}
Let
\begin{equation*}
(\mathcal{M}, \phi) = (A_1, \phi_1) \ast (A_2, \phi_2)
\end{equation*}
be the von Neumann algebra free product of factors of type ${\rm I}$ with respect to faithful states, at least one of which is nontracial. Then $\mathcal{M}$ is a full factor of type ${\rm III}$ and $\phi$ is an almost periodic faithful state whose centralizer is isomorphic to the type ${\rm II_1}$ factor $L(\F_{\infty})$. The point spectrum of the modular operator $\Delta_{\phi}$ of $\phi$, is equal to the subgroup of $\R^*_+$ generated by the union of the point spectra of $\Delta_{\phi_1}$ and of $\Delta_{\phi_2}$. Thus in Connes' classification, $\mathcal{M}$ is always a factor of type ${\rm III_{\lambda}}$, with $0 < \lambda \leq 1$.
\end{theo}
The fact that $\phi$ is an almost periodic faithful state \cite{connes74} is an easy consequence of basic results on free products (see Section $4$ in \cite{dykema96} for further details). Of course, the fact that the centralizer of $\phi$ is isomorphic to the type ${\rm II_1}$ factor $L(\F_{\infty})$ is the most difficult part of the theorem. To prove this, Dykema uses sophisticated algebraic techniques on free products that he developed also in \cite{dykema94} and \cite{dykema93}. Finally, the fact that $\mathcal{M}$ is of type ${\rm III}$ follows from results of Section $4$ in \cite{dykema96}. 

However, natural questions were left unanswered. Dykema asked in Question $9.1$ \cite{dykema96}, whether the type ${\rm III_{\lambda}}$ factors that are obtainable by taking various free products of finite dimensional or hyperfinite algebras are isomorphic to each other, and more precisely whether they are isomorphic to the factor of R\u{a}dulescu \cite{radulescu1993},
\begin{equation*}
(L(\Z), \tau_{\Z}) \ast (M_2(\C), \omega_{\lambda})
\end{equation*}
where $\omega_{\lambda}(p_{ij}) = \delta_{ij}\lambda^j/(\lambda + 1)$ for $i, j \in \{0, 1\}$. Furthermore, he asked in Question $9.3$ \cite{dykema96}, whether the full factors of type ${\rm III_1}$ having the same $\Sd$ invariant that are obtainable by taking free products of various finite dimensional or hyperfinite algebras are isomorphic to each other. We will see that we partially give positive answers to these questions.

In \cite{shlya97}, Shlyakhtenko introduced a new class of full factors of type ${\rm III}$. His idea is to give  a version of the CAR\footnote{Canonical Anticommutation Relations} functor and of the associated quasi-free states in the framework of Voiculescu's free probability theory \cite{voiculescu92}. In Section $2$, we recall his construction. But roughly speaking, we can say that to each real Hilbert space $H_{\R}$ and to each orthogonal representation $(U_t)$ of $\R$ on $H_{\R}$, he associated a factor $\Gamma(H_\R, U_t)''$ called the \emph{free Araki-Woods factor}. He proved that $\Gamma(H_\R, U_t)''$ is a type ${\rm III}$ factor except if $U_t = \id$ for all $t \in \R$. The restriction to $\Gamma(H_\R, U_t)''$ of the vacuum state denoted by $\varphi_U$ and called the \emph{free quasi-free state}, is faithful. Moreover, he proved that $\varphi_U$ is an almost periodic state iff the orthogonal representation $(U_t)$ is almost periodic. Recall in this respect the following definition:

\begin{df} [Connes, \cite{connes74}]
\emph{Let $M$ be a von Neumann algebra with separable predual which has almost periodic weights. The \emph{$\Sd$ invariant} of $M$ is defined as the intersection over all the almost periodic, faithful, normal, semifinite weights $\varphi$ of the point spectra of the modular operators $\Delta_{\varphi}$.}
\end{df}
Connes proved that for a factor of type ${\rm III}$, $\Sd(M)$ is a countable subgroup of $\R^*_+$ \cite{connes74}. 
In the almost periodic case, using a powerful tool called the \emph{matricial model}, Shlyakhtenko obtains this remarkable result:
\begin{theo}[Shlyakhtenko, \cite{{shlya2004}, {shlya97}}]\label{pparaki}
Let $(U_t)$ be a nontrivial almost periodic orthogonal representation of $\R$ on the real Hilbert space $H_{\R}$ with $\dim H_{\R} \geq 2$. Let $A$ be the infinitesimal generator of $(U_t)$ on $H_{\C}$, the  complexified Hilbert space of $H_{\R}$. Denote $M = \Gamma(H_{\R}, U_t)''$. Let $\Gamma \subset \R^*_+$ be the subgroup generated by the point spectrum of $A$. Then, $M$ only depends on $\Gamma$ up to state-preserving isomorphisms. 

Conversely, the group $\Gamma$ coincides with the $\Sd$ invariant of the factor $M$. Consequently, $\Sd$ completely classifies the almost periodic free Araki-Woods factors. Moreover, the centralizer of the free quasi-free state $\varphi_U$ is isomorphic to the type ${\rm II_1}$ factor $L(\F_{\infty})$.
\end{theo}
He proved also that the (unique) free Araki-Woods factor of type ${\rm III_{\lambda}}$, denoted by $(T_{\lambda}, \varphi_{\lambda})$, is isomorphic to the factor of R\u{a}dulescu \cite{radulescu1993} and {\textquotedblleft freely absorbs\textquotedblright} $L(\F_{\infty})$. Since the free Araki-Woods factors satisfy free absorption properties, Shlyakhtenko asked whether the free products of matrix algebras $(A_1, \phi_1) \ast (A_2, \phi_2)$ are stable by taking free products with $L(\Z)$, in other words whether they are free Araki-Woods factors. 

We give in this paper a positive answer to the question of Shlyakhtenko for certain free products of matrix algebras and other von Neumann algebras. Thanks to Theorem $6.6$ of \cite{shlya97}, it partially gives positive answers to Questions $9.1$ and $9.3$ of Dykema \cite{dykema96}. Before giving the results, we must introduce a few notations. For an almost periodic state $\phi$, we denote by $\Sd(\phi)$ the subgroup of $\R^*_+$ generated by the point spectrum of the modular operator $\Delta_{\phi}$. On $B(\ell^2(\N))$, we denote by $\psi_{\lambda}$ the state given by $\psi_{\lambda}(e_{ij}) = \delta_{ij}\lambda^j(1 - \lambda)$ for $i, j \in \N$. For $\beta \in ]0, 1[$, we denote by $(\C^2, \tau_{\beta})$ the algebra generated by a projection $q$ with $\tau_{\beta}(q) = \beta$. The hyperfinite type ${\rm II_1}$ factor together with its trace is denoted by $(\mathcal{R}, \tau)$. At last, we denote by $(T_{\Gamma}, \varphi_{\Gamma})$ the unique (up to state-preserving isomorphism) almost periodic free Araki-Woods factor whose $\Sd$ invariant is exactly $\Gamma$ (see Notation $\ref{araki}$ for further details). 

\begin{df}
\emph{Let $\rho : (B, \phi_B) \hookrightarrow (A, \phi_A)$ be an embedding of von Neumann algebras. We shall say that $\rho$ is \emph{modular} if it is state-preserving and if $\rho(B)$ is globally invariant under the modular group $(\sigma_t^{\phi_A})$.}  
\end{df}
Thanks to a result of Takesaki (see for example Theorem ${\rm IX}.4.2$ in \cite{takesakiII}), $\rho : (B, \phi_B) \hookrightarrow (A, \phi_A)$ is modular if and only if there exists a state-preserving conditional expectation from $A$ onto $\rho(B)$. Note that in this case, such a conditional expectation is unique. Our main results are given in the following theorem.

\begin{theo}\label{results}
Let $(A_i, \phi_i)$, $i = 1, 2$, be two von Neumann algebras endowed with a faithful, normal, almost periodic state $\phi_i$, such that for $i = 1, 2$
\begin{equation*}
(A_i, \phi_i) \ast (L(\Z), \tau_\Z) \cong (T_{\Sd(\phi_i)}, \varphi_{\Sd(\phi_i)}).
\end{equation*}
Let $\Gamma$ be the subgroup of $\R^*_+$ generated by $\Sd(\phi_1)$ and $\Sd(\phi_2)$.  Assume that for some $\lambda, \beta \in ]0, 1[$, there exist modular embeddings
\begin{eqnarray*}
(M_2(\C), \omega_\lambda) & \hookrightarrow & (A_1, \phi_1) \\
(\C^2, \tau_\beta) & \hookrightarrow & (A_2, \phi_2),
\end{eqnarray*}
such that $\lambda/(\lambda + 1) \leq \min\{\beta, 1 - \beta\}$. Then
\begin{equation*}
(T_{\Gamma}, \varphi_\Gamma) \cong (A_1, \phi_1) \ast (A_2, \phi_2).
\end{equation*}
\end{theo}

In particular, for any  $\lambda, \mu \in ]0, 1[$, $(M_2(\C), \omega_{\lambda}) \ast 
(M_2(\C), \omega_{\mu})$, $(M_2(\C), \omega_{\lambda}) \ast 
(\mathcal R, \tau)$ and $(B(\ell^2(\N)), \psi_{\lambda}) \ast 
(\mathcal R, \tau)$ are free Araki-Woods factors.

The paper is organized as follows. Section $\ref{bac}$ is devoted to a few reminders on free products and free Araki-Woods factors. In Section $\ref{study}$, we show that $\displaystyle{(M_2(\C), \omega_{\lambda}) \ast (\C^2, \tau_\beta)}$ is isomorphic in a state-preserving way to $(T_{\lambda}, \varphi_{\lambda})$, whenever $\lambda/(\lambda + 1) \leq \min\{\beta, 1 - \beta\}$. In Section $\ref{mainthe}$, we prove Theorem $\ref{results}$ using the {\textquotedblleft machinery\textquotedblright} of amalgamated free products. Finally, Section $\ref{rema}$ is a remark.

\section{Conventions and Preliminary Background}\label{bac}

Throughout this paper, we will be working with free products of von Neumann algebras with respect to states. For the convenience of the reader, it is useful to remind the following notation:
\begin{nota}
\emph{If $(M, \varphi)$ and $(N, \psi)$ are von Neumann algebras endowed with states $\varphi$ and $\psi$, the notation $(M, \varphi) \cong (N, \psi)$ means that there exists a $\ast$-isomorphism $\alpha : M \to N$ such that $\psi\alpha = \varphi$.}
\end{nota}
We remind this well-known proposition concerning free products of von Neumann algebras with respect to states.
\begin{prop}[\cite{voiculescu92}]\label{libre}
Let $(M_i, \varphi_i)$ be a family of von Neumann algebras endowed with faithful normal states. Then, there exists, up to state-preserving isomorphism, a unique von Neumann algebra $(M, \varphi)$ endowed with a faithful normal state $\varphi$ such that
\begin{enumerate}
\item [$-$] $(M_i, \varphi_i)$ embeds into $(M, \varphi)$ in a state-preserving way,
\item [$-$] $M$ is generated by the family of subalgebras $(M_i)$ which is a free family in $(M, \varphi)$.
\end{enumerate}
The \emph{free product} of $(M_i, \varphi_i)$ is denoted by $(M, \varphi) = \displaystyle{\mathop{\ast}_{i \in I}} (M_i, \varphi_i)$.
\end{prop}
 
\begin{nota}[\cite{dykema93}]
\emph{For von Neumann algebras $A$ and $B$, with states $\varphi_A$ and $\varphi_B$, the von Neumann algebra}
$$\mathop{A}_{\alpha}^{p} \oplus \mathop{B}_{\beta}^{q}$$
\emph{where $\alpha, \beta \geq 0$ and $\alpha + \beta = 1$, will denote the algebra $A \oplus B$ whose associated state is $\varphi(a, b) = \alpha\varphi_A(a) + \beta\varphi_B(b)$. Moreover, $p \in A$ and $q \in B$ are projections corresponding to the identity elements of $A$ and $B$.}
\end{nota}

Now, we want to remind the construction of the free Araki-Woods factors \cite{shlya97}. Let $H_{\R}$ be a real Hilbert space and let $(U_t)$ be an orthogonal representation of $\R$ on $H_{\R}$. Let $H = H_{\R} \otimes_{\R} \C$ be the complexified Hilbert space. If $A$ is the infinitesimal generator of $(U_t)$ on $H$, we remind that $j : H_{\R} \to H$ defined by $j(\zeta) = (\frac{2}{A^{-1} + 1})^{1/2}\zeta$ is an isometric embedding of $H_{\R}$ into $H$. Let $K_{\R} = j(H_{\R})$. Introduce the \emph{full Fock space} of $H$:
$$\mathcal{F}(H) =\C\Omega \oplus \bigoplus_{n = 1}^{\infty} H^{\otimes n}.$$ 
The unit vector $\Omega$ is called \emph{vacuum vector}. For any $\xi \in H$, we have the left creation operator
$$l(\xi) : \mathcal{F}(H) \to \mathcal{F}(H) : \left\{ 
{\begin{array}{l} l(\xi)\Omega = \xi, \\ l(\xi)(\xi_1 \otimes \cdots \otimes \xi_n) = \xi \otimes \xi_1 \otimes \cdots \otimes \xi_n.
\end{array}} \right.
$$
For any $\xi \in H$, we denote by $s(\xi)$ the real part of $l(\xi)$ given by
$$s(\xi) = \frac{l(\xi) + l(\xi)^*}{2}.$$
The crucial result of Voiculescu \cite{voiculescu92} claims that the distribution of the operator $s(\xi)$ with respect to the vacuum vector state $\varphi(x) = \langle x\Omega, \Omega\rangle$ is the semicircular law of Wigner supported on the interval $[-\|\xi\|, \|\xi\|]$. 

\begin{df}[Shlyakhtenko, \cite{shlya97}]
\emph{Let $(U_t)$ be an orthogonal representation of $\R$ on the real Hilbert space $H_{\R}$ $(\dim H_{\R} \geq 2)$. The \emph{free Araki-Woods factor} associated with $H_\R$ and $(U_t)$, denoted by $\Gamma(H_{\R}, U_t)''$, is defined by}
$$\Gamma(H_{\R}, U_t)'' = \{s(\xi), \xi \in K_{\R}\}''.$$
\emph{The vector state $\varphi_{U}(x) = \langle x\Omega, \Omega\rangle$ is called the {\it free quasi-free state}.}
\end{df}
As we said previously, the free Araki-Woods factors provide many new examples of full factors of type {\rm III} \cite{{barnett95}, {connes73}, {shlya2004}}. We can summarize the general properties of free Araki-Woods factors in the following theorem (see also \cite{vaes2004}):

\begin{theo}[Shlyakhtenko, \cite{{shlya2004}, {shlya99}, {shlya98}, {shlya97}}]
Let $(U_t)$ be an orthogonal representation of $\R$ on the real Hilbert space $H_{\R}$ with $\dim H_{\R} \geq 2$. Denote $M = \Gamma(H_{\R}, U_t)''$.
\begin{enumerate}
\item $M$ is a full factor.
\item $M$ is of type ${\rm II_1}$ iff $U_t = id$ for every $t \in \R$.
\item $M$ is of type ${\rm III_{\lambda}}$ $(0 < \lambda < 1)$ iff $(U_t)$ is periodic of period $\frac{2\pi}{|\log \lambda|}$.
\item $M$ is of type ${\rm III_1}$ in the other cases.
\item The factor $M$ has almost periodic states iff $(U_t)$ is almost periodic.
\end{enumerate}
\end{theo}
Let $H_{\R} = \R^2$ and $0 < \lambda < 1$. Let
\begin{equation}\label{periodic}
U_t = \begin{pmatrix}
\cos(t\log \lambda) & - \sin(t\log \lambda) \\
\sin(t\log \lambda) & \cos(t\log \lambda)
\end{pmatrix}.
\end{equation}

\begin{nota}[\cite{shlya97}]\label{Tlambda}
\emph{Denote $(T_{\lambda}, \varphi_{\lambda}) := \Gamma(H_{\R}, U_t)''$ where $H_{\R} = \R^2$ and $(U_t)$ is given by Equation $(\ref{periodic})$.}
\end{nota}
Using a powerful tool called the \emph{matricial model}, Shlyakhtenko was able to prove the following isomorphism
$$(T_{\lambda}, \varphi_{\lambda}) \cong (B(\ell^2(\N)), \psi_{\lambda}) \ast (L^{\infty}[-1, 1], \mu),$$
where $\psi_{\lambda}(e_{ij}) = \delta_{ij}\lambda^j(1 - \lambda)$, $i, j \in \N$, and $\mu$ is a nonatomic measure on $[-1, 1]$. He also proved that $(T_{\lambda}, \varphi_{\lambda})$ is isomorphic to the factor of type ${\rm III_{\lambda}}$ introduced by R\u{a}dulescu in \cite{radulescu1993}. Namely,
$$(T_{\lambda}, \varphi_{\lambda}) \cong (M_2(\C), \omega_{\lambda}) \ast (L^{\infty}[-1, 1], \mu),$$
where $\omega_{\lambda}(p_{ij}) = \delta_{ij}\lambda^j/(\lambda + 1)$, $i, j \in \{0, 1\}$, and $\mu$ a nonatomic measure on $[-1, 1]$. Moreover, he showed that $(T_{\lambda}, \varphi_{\lambda})$ has a good behaviour when it is compressed by a {\textquotedblleft right\textquotedblright} projection. More precisely, denote $(C, \psi) := (B(\ell^2(\N)), \psi_{\lambda}) \ast (L^{\infty}[-1, 1], \mu)$ and $(D, \omega) := (M_2(\C), \omega_{\lambda}) \ast (L^{\infty}[-1, 1], \mu)$. The following proposition is an easy consequence of proofs of Theorems $5.4$ and $6.7$ of \cite{shlya97}. It will be useful in Section $\ref{study}$.

\begin{prop}\label{compression}
Let $(C, \psi)$, $(D, \omega)$ defined as above and $e_{00} \in B(\ell^2(\N)) \subset C$, $p_{00}$,  $p_{11} \in M_2(\C) \subset D$. Then
\begin{eqnarray*}
(T_{\lambda}, \varphi_{\lambda}) & \cong & (e_{00}Ce_{00}, \frac{1}{\psi(e_{00})}\psi) \\
& \cong & (p_{00}Dp_{00}, \frac{1}{\omega(p_{00})}\omega) \\
& \cong & (p_{11}Dp_{11}, \frac{1}{\omega(p_{11})}\omega).
\end{eqnarray*}
\end{prop}

When the representation $(U_t)$ is assumed to be almost periodic, we have seen (Theorem $\ref{pparaki}$) that $\Gamma \subset \R^*_+$, the subgroup generated by the point spectrum of $A$, completely classifies the free Araki-Woods factor $\Gamma(U_t, H_{\R})''$. 
\begin{nota}\label{araki}
\emph{For any nontrivial countable subgroup $\Gamma \subset \R^*_+$, we shall denote by $(T_{\Gamma}, \varphi_{\Gamma})$ the unique (up to state-preserving isomorphism) almost periodic free Araki-Woods factor whose $\Sd$ invariant is exactly $\Gamma$. Of course, $\varphi_{\Gamma}$ is its free quasi-free state. If $\Gamma = \lambda^{\Z}$ for $\lambda \in ]0, 1[$, then $(T_{\Gamma}, \varphi_{\Gamma})$ is of type ${\rm III_{\lambda}}$; in this case, it will be simply denoted by $(T_{\lambda}, \varphi_{\lambda})$ \cite{shlya97}, as in Notation $\ref{Tlambda}$. Theorem $6.4$ in \cite{shlya97} gives the following formula:}
$$(T_{\Gamma}, \varphi_{\Gamma}) \cong \displaystyle{\mathop{\ast} _{\gamma \in \Gamma} (T_{\gamma}, \varphi_{\gamma})}.$$
\end{nota}

\section{Study of $\displaystyle{(M_2(\C), \omega_{\lambda}) \ast \left(\mathop{\C}_{\beta} \oplus \mathop{\C}_{1 - \beta}\right)}$ for $\lambda/(\lambda + 1) \leq \min\{\beta, 1 - \beta\}$}\label{study}

For any $\beta \in ]0, 1[$, the von Neumann algebra $\displaystyle{\mathop{\C}_{\beta} \oplus \mathop{\C}_{1 - \beta}}$ is simply denoted by $(\C^2, \tau_{\beta})$. Let $\lambda \in ]0, 1[$ and denote $\alpha = \lambda/(\lambda +1)$. We remind that the faithful state $\omega_{\lambda}$ on $M_2(\C)$ is defined as follows: $\omega_{\lambda}(p_{ij}) = \delta_{ij}\lambda^j/(\lambda +1)$, for $i, j \in \{0, 1\}$. The aim of this section is to prove the following theorem:
\begin{theo}\label{premierresultat}
If $\alpha = \lambda/(\lambda +1) \leq \min\{\beta, 1 - \beta\}$, then
\begin{equation}\label{prequation}
(M_2(\C), \omega_{\lambda}) \ast (\C^2, \tau_{\beta}) \cong (T_{\lambda}, \varphi_{\lambda}).
\end{equation}
\end{theo}
\begin{nota}
\emph{The von Neumann algebra of the left-hand side of $(\ref{prequation})$ together with its free product state will be denoted by $(M, \omega)$.}
\end{nota}
To prove Theorem $\ref{premierresultat}$, we will need the following result due to Voiculescu \cite{voiculescu87} (see also \cite{{dykema94}, {dykema93}}) which gives a precise picture of the von Neumann algebra generated by two projections $p$ and $q$ free with respect to a faithful trace. More precisely,
\begin{theo}[Voiculescu, \cite{voiculescu87}]\label{deuxprojs}
Let $0 < \alpha \leq \min\{\beta, 1 - \beta\} < 1$. Then
\begin{equation}\label{deuxproj}
\left( \mathop{\C}_{\alpha}^{p} \oplus \mathop{\C}_{1 - \alpha}^{1 - p}\right) \ast \left( \mathop{\C}_{\beta}^{q} \oplus \mathop{\C}_{1 - \beta}^{1 - q} \right) \cong \mathop{\C}_{\beta - \alpha}^{(1 - p) \wedge q} \oplus \mathop{\underbrace{\left(L^{\infty}\left(\left[0, \frac{\pi}{2}\right], \nu \right) \otimes M_2\right)}}_{2\alpha} \oplus \mathop{\C}_{1 - \alpha - \beta}^{(1 - p) \wedge (1 - q)},
\end{equation}
where $\nu$ is a probability measure without atoms on $[0, \pi/2]$, and $L^{\infty}([0, \pi/2], \nu)$ has trace given by integration against $\nu$. In the picture of the right-hand side of $(\ref{deuxproj})$, we have
$$p = 0 \oplus 
\begin{pmatrix}
1 & 0 \\ 
0 & 0
\end{pmatrix} \oplus 0,$$
$$q = 1 \oplus 
\begin{pmatrix}
\cos^2\theta & \cos\theta\sin\theta \\ 
\cos\theta\sin\theta & \sin^2\theta
\end{pmatrix}
\oplus 0,$$
where $\theta \in [0, \pi/2]$.
\end{theo}
\begin{rem}\label{generateurs}
\emph{The von Neumann algebra of the right-hand side of $(\ref{deuxproj})$ together with its trace will be denoted by $(N, \tau)$.
With the previous notations, $(N, \tau)$ embeds into $(M, \omega)$ in a state-preserving way. We shall assume that $N \subset M$, and $p = p_{11}$. Consequently, $p$ is the {\textquotedblleft smallest\textquotedblright} projection with respect to the state $\omega$, i.e. $\omega(p) = \alpha = \lambda/(\lambda + 1)$. As in \cite{dykema94}, let
$$x = 0 \oplus 
\begin{pmatrix}
0 & 0 \\
1 & 0
\end{pmatrix}
\oplus 0 = \pol((1 - p)qp),$$
where {\textquotedblleft pol\textquotedblright} means {\textquotedblleft polar part of\textquotedblright}. Then $x$ is a partial isometry from $p$ into $1 - p$, i.e., $x^*x = p$ and $xx^* \leq 1 - p$. Let $z$ be the projection
$$z = 1 \oplus 
\begin{pmatrix}
0 & 0 \\
0 & 1
\end{pmatrix}
\oplus 0.$$
Note that $z \leq 1 - p$. We see that $N$ is the von Neumann generated by $pqp$ together with $x$ and $z$. We shall denote $N = W^*(pqp, x, z)$.
}
\end{rem}
\begin{df}[\cite{dykema94}]
\emph{Let $(S_\iota)_{\iota \in I}$ be a family of subsets of a unital algebra $A \owns 1$. A nontrivial \emph{traveling product} in $(S_\iota)_{\iota \in I}$ is a product $a_1 \cdots a_n$ such that $a_j \in S_{\iota_j}$ $(1 \leq j \leq n)$ and $\iota_1 \neq \iota_2 \neq \cdots \neq \iota_{n - 1} \neq \iota_n$. The \emph{trivial traveling product} is the identity element $1$. The set of all traveling products in $(S_{\iota})_{\iota \in I}$, including the trivial one is denoted by $\Lambda((S_{\iota})_{\iota \in I})$. If $|I| = 2$, we will call traveling products {\it alternating products}.}
\end{df}
We are now ready to prove the following proposition; it gives a precise picture of the compression of the von Neumann algebra $(M, \omega)$ by the projection $p$. The proof is based on algebraic techniques developed in \cite{{dykema96}, {dykema94}, {dykema93}}, and techniques of computation of $\ast-$distributions developed in \cite{shlya97} and \cite{voiculescu90}.
\begin{prop}\label{comp}
Let $(M, \omega) = (M_2(\C), \omega_{\lambda}) \ast (\C^2, \tau_{\beta})$ and $p = p_{11} \in M_2(\C)$. Assume as in Theorem $\ref{deuxprojs}$, that $\alpha = \lambda/(\lambda + 1) \leq \min\{\beta, 1 - \beta\}$. Then
\begin{equation*}
\left(pMp, \frac{1}{\omega(p)}\omega\right)  \cong  L(\Z) \ast \left((\C^2, \tau_{\delta}) \otimes (B(\ell^2(\N)), \psi_{\lambda})\right),
\end{equation*}
where $(\C^2, \tau_{\delta}) = \displaystyle{\mathop{\C}_{\delta} \oplus \mathop{\C}_{1 - \delta}}$ with $\delta = \frac{1 - \frac{\beta}{1 - \alpha}}{1 - \lambda}$ and $\psi_{\lambda}(e_{ij}) = \delta_{ij}\lambda^j(1 - \lambda)$, for $i$, $j \in \N$.
\end{prop}
\begin{proof}
Let $(M, \omega) = (M_2, \omega_{\lambda}) \ast (\C^2, \tau_{\beta})$. Let $p$ and $q$ be the projections in $M$ such that $N = W^*(p, q)$ as in Theorem $\ref{deuxprojs}$; $p$ and $q$ are free in $M$ with respect to $\omega$ and $\omega(p) = \alpha = \lambda/(\lambda + 1)$, $\omega(q) = \beta$. Let $x$ and $z$ as in Remark $\ref{generateurs}$. We know that $N = W^*(pqp, x, z)$. Denote by $u = p_{01} \in M_2(\C)$ the partial isometry from $p$ to $1 - p$, i.e. $u^*u = p$ and $uu^* = 1 - p$. Then, thanks to Lemma 5.3 from \cite{voiculescu92}
\begin{equation*}
pMp = W^*(pqp, u^*x, u^*zu).
\end{equation*}
Denote $v = u^*x$ and $P = u^*zu$. Since $v^*v = p$ and $vv^* \leq p$, $v$ is an isometry in $pMp$. Moreover, since $Pv = u^*zuu^*x = u^*zx = u^*x = v$, we get $vv^* \leq P$. Denote $\omega_p = \frac{1}{\omega(p)}\omega$ the canonical state on $pMp$. First, we are going to compute the $\ast-$distributions of the elements $v$ and $vP$ in $pMp$ with respect to $\omega_p$.
\begin{lem}\label{distribution}
Let $\gamma = \beta(\lambda + 1) = \beta/(1 - \alpha)$. For any $k, l \in \N$,
\begin{eqnarray}\label{distribution1}
\omega_p(v^k(v^*)^l) & = & \delta_{kl}\lambda^k, \\ 
\omega_p(v^kP(v^*)^l) & = & \delta_{kl}\lambda^k\gamma. \label{distribution2}
\end{eqnarray}
\end{lem}
\begin{proof}[Proof of Lemma $\ref{distribution}$]
{\bf Step (0)}. First, we review the {\textquotedblleft algebraic trick\textquotedblright} of Dykema \cite{dykema94}. Denote $a = p - \omega(p)$ and $b = q - \omega(q)$; we have $N = \overline{\mbox{span}}^w\Lambda(\{a\}, \{b\})$. Let $w \in N$ such that $\omega(w) = \omega(pw) = 0$. By Kaplansky Density Theorem, $w$ is the s.o.-limit of a bounded sequence in $\mbox{span}\:\Lambda(\{a\}, \{b\})$. Note that since $a$ and $b$ are free and $\omega(a) = \omega(b) = 0$, if $y \in \mbox{span}\:\Lambda(\{a\}, \{b\})$, then $\omega(y)$ is equal to the coefficient of $1$ in $y$. Since $\omega(w) = 0$, we may choose that approximating sequence in $\mbox{span}\: \left(\Lambda(\{a\}, \{b\}) \backslash \{1\}\right)$. Moreover, since $\omega(pw) = 0$, we may also assume that each coefficient of $a$ be zero, i.e. we have a bounded approximating sequence for $w$ of elements of $\mbox{span} \: (\Lambda(\{a\}, \{b\}) \backslash \{1, a\})$.

{\bf Step (1)}. We prove now Equation $(\ref{distribution1})$. Assume $k \geq 1$ and $l = 0$, then $v^k = (u^*x)^k$ is a nontrivial alternating product in $\{u^*\}$ and $\{x\}$. Since $\omega(x) = \omega(px) = 0$, $x$ is a s.o.-limit of a bounded sequence in $\mbox{span} \left(\Lambda(\{a\}, \{b\}) \backslash \{1, a\}\right)$. So to show that $\omega_p(v^k) = 0$, it suffices to show that if $s$ is a nontrivial alternating product in $\mbox{span} \left(\Lambda(\{a\}, \{b\}) \backslash \{1, a\}\right)$ and $\{u^*\}$ then $\omega(s) = 0$. But since $u^*a = -\alpha u^*$ and $au^* = (1 - \alpha)u^*$, regrouping gives a nontrivial alternating product in $\{a, u^*\}$ and $\{b\}$, hence by freeness $\omega(s) = 0$. We get also immediatly $\omega_p((v^*)^l) = 0$. Assume at last $k \geq 1$ and $l \geq 1$, then $v^k(v^*)^l = (u^*x)^{k - 1}u^*xx^*u(x^*u)^{l - 1}$. Let $y = xx^* - \alpha1 + \lambda a$. Since $\omega(xx^*) = \omega(x^*x) =\alpha$ and $py = 0$, $\omega(y) = \omega(py) = 0$, hence $y$ is a s.o.-limit of a bounded sequence in $\mbox{span}\: \left(\Lambda(\{a\}, \{b\}) \backslash \{1, a\}\right)$. Replacing in $(u^*x)^{k - 1}u^*xx^*u(x^*u)^{l - 1}$ the term $xx^*$ by $y + \alpha1 - \lambda a$, and since $u^*au = -\alpha p$, we have
\begin{equation}\label{decomposition1}
\omega_p(v^k(v^*)^l) = \omega_p((u^*x)^{k - 1}u^*yu(x^*u)^{l - 1}) + \lambda\omega_p((u^*x)^{k - 1}p(x^*u)^{l - 1}).
\end{equation}
To prove $\omega_p((u^*x)^{k - 1}u^*yu(x^*u)^{l - 1}) = 0$, it suffices to show that if $r$ is a nontrivial alternating product in $\mbox{span} \left(\Lambda(\{a\}, \{b\}) \backslash \{1, a\}\right)$ and $\{u^*, u\}$ then $\omega(s) = 0$. But for the same reasons as above, regrouping gives a nontrivial alternating product in $\{a, u^*, u\}$ and $\{b\}$, hence by freeness $\omega(r) = 0$. If in Equation $(\ref{distribution1})$, $k \neq l$, then applying Equation $(\ref{decomposition1})$ several times we eventually get $\omega_p(u^*x \cdots u^*x)$ or $\omega_p(x^*u \cdots x^*u)$, both of which are zero. If $k = l$, then we eventually get $\lambda^k\omega_p(p) = \lambda^k$. Thus Equation $(\ref{distribution1})$ holds.

{\bf Step (2)}. We prove at last Equation $(\ref{distribution2})$. Since $\omega(u^*zu) = \lambda\omega(z)$, and $\gamma =  \beta(\lambda + 1) = \beta/(1 - \alpha)$, we get $\omega_{p}(P) = \gamma$. Assume $k$, $l \geq 0$, then $v^kP(v^*)^l = (u^*x)^ku^*zu(x^*u)^l$. Since $\omega(z) = \beta$ and $pz = 0$, $y = z - \beta1 + \gamma a$ satisfies $\omega(y) = \omega(py) = 0$.  Consequently, $y$ is a s.o.-limit of a bounded sequence in $\mbox{span}(\Lambda(\{a\}, \{b\} \backslash \{1, a\})$. Replacing in the product $(u^*x)^ku^*zu(x^*u)^l$ the term $z$ by $y + \beta1 - \gamma a$, and since $u^*au = -\alpha p$, we have
\begin{equation*}
\omega_p(v^kP(v^*)^l) = \omega_p((u^*x)^{k}u^*yu(x^*u)^{l}) + \gamma\omega_p((u^*x)^{k}p(x^*u)^{l}).
\end{equation*}
For the same reasons, $\omega_p((u^*x)^{k}u^*yu(x^*u)^{l}) = 0$ and  $\omega_p(v^kP(v^*)^l) = \gamma\omega_p(v^{k}(v^*)^{l})$. Thus Equation $(\ref{distribution2})$ holds.
\end{proof}

\begin{lem}\label{freeness}
In $pMp$, the von Neumann subalgebras $W^*(pqp)$ and $W^*(v, P)$ are $\ast$-free with respect to $\omega_p$.
\end{lem}
\begin{proof}[Proof of Lemma $\ref{freeness}$]
Lemma $4.2$ in \cite{shlya97} inspired us to prove Lemma $\ref{freeness}$. It is slightly more complicated because here, in some sense, the assumptions are weaker and we must additionally deal with the projection $P$. To overcome these difficulties, we will use the {\textquotedblleft algebraic trick\textquotedblright} \cite{dykema94} mentioned above.

Let $B = W^*(pqp)$ be the von Neumann subalgebra of $pMp$ generated by $pqp$ and $C = W^*(v, P)$ be the von Neumann subalgebra of $pMp$ generated by $v$ and $P$. Let $g_k = (pqp)^k - \omega_p((pqp)^k)p$ for $k \geq 1$. Let $W_{kl} = v^k(v^*)^l - \delta_{kl}\lambda^kp$, $W'_{rs} = v^rP(v^*)^s - \gamma\delta_{rs}\lambda^rp$ for $k, l, r, s \in \N$, $k + l > 0$. Since 
\begin{eqnarray*}
B & = & \overline{\mbox{span}}^w \{p, g_k \:|\: k \geq 1\}, \\
C & = &  \overline{\mbox{span}}^w\{p, W_{kl}, W'_{rs} \:|\: k, l, r, s \in \N, k + l > 0 \}, 
\end{eqnarray*}
it follows that to check freeness of $B$ and $C$, we must show that
\begin{equation}\label{but}
\omega_p(\mathop{\underbrace{b_0w_1b_1 \cdots w_nb_n}}_W) = 0
\end{equation}
where
\begin{eqnarray}\label{word1}
b_j & = & g_{m_j}, \\ \label{word2}
w_j & = & W_{k_jl_j} \\ \label{word3}
\mbox{ or } w_j & = & W'_{r_js_j},
\end{eqnarray}
with $k_j, l_j, m_j, r_j, s_j \in \N, k_j + l_j > 0, m_j > 0$ for all $j$, except possibly $b_0$ and/or $b_n$ are equal to $1$. We shall prove Equation (\ref{but}) under a weaker assumption, which is that $w_j$ is also allowed to be 
\begin{equation}\label{reduction}
(u^*x)^{s_j}u^*yu(x^*u)^{t_j},
\end{equation}
$s_j, t_j \geq 0$ and $y = xx^* - \alpha1 + \lambda a$ or $y = z - \beta1 + \gamma a$ as in proof of Lemma \ref{distribution}. 

We will denote by $W = b_0w_1b_1 \cdots w_nb_n$ such a word with $w_j$ as in Equation $(\ref{word2})$ or $(\ref{word3})$. Let $w_j$ be as in Equation $(\ref{word2})$ with both $k_j$ and $l_j$ nonzero and let $y_1 = xx^* - \alpha1 + \lambda a$. We will replace this $w_j$ by
\begin{eqnarray*}
w_j & = & ((u^*x)^{k_j - 1}u^*y_1u(x^*u)^{l_j - 1}) \\ \nonumber
 & & + \: (\lambda(u^*x)^{k_j - 1}(x^*u)^{l_j - 1} - \delta_{k_jl_j}\lambda^{k_j}) \\ \nonumber
 & = & A_j + B_j. \nonumber
\end{eqnarray*}
Let now $w_j$ be as in Equation $(\ref{word3})$ with both $r_j$ and $s_j$ nonzero and let $y_2 = z - \beta1 +  \gamma a$. We will replace this $w_j$ by
\begin{eqnarray*}
w_j & = & \left((u^*x)^{r_j}u^*y_2u(x^*u)^{s_j}\right) \\ \nonumber
 & & + \: \gamma\left(\lambda(u^*x)^{r_j}(x^*u)^{s_j} - \delta_{r_js_j}\lambda^{r_j}\right) \\ \nonumber
 & = & \left((u^*x)^{r_j}u^*y_2u(x^*u)^{s_j}\right) \\ \nonumber
 & & + \: \gamma\left((u^*x)^{r_j - 1}u^*y_1u(x^*u)^{s_j - 1}\right) \\ \nonumber
 & & + \: \gamma\left(\lambda(u^*x)^{r_j - 1}(x^*u)^{s_j - 1} - \delta_{r_js_j}\lambda^{r_j}\right) \\ \nonumber
 & = & A'_j + A''_j + C_j.\nonumber 
\end{eqnarray*}

After such replacements are done, $w$ can be rewritten as a sum of terms, in which some $w_j$ are replaced by $A_j$'s, $A'_j$'s, $A''_j$'s, some by $B_j$'s and some by $C_j$'s. Consider the terms where all replacements are replacements by $A_j$'s, $A'_j$'s, $A''_j$'s. These terms can be written as alternating products  in $\Omega = \{x, x^*, y_1, y_2, g_k, xg_k, g_kx^*, xg_kx^*\}$ and $\{u, u^*\}$. But each element $h \in \Omega$ satisfies $\omega(h) = \omega(ph) = 0$, hence $h$ is a s.o.-limit of a bounded sequence in $\mbox{span} \left(\Lambda(\{a\}, \{b\}) \backslash \{1, a\}\right)$. We use the same argument as before. To prove that $\omega_p$ is zero on such terms, it suffices to show that $\omega$ is zero on a nontrivial alternating product in $\mbox{span} \left(\Lambda(\{a\}, \{b\}) \backslash \{1, a\}\right)$ and $\{u, u^*\}$. But regrouping gives a nontrivial alternating product in $\{a, u, u^*\}$ and $\{b\}$. So, by freeness $\omega$ is zero on such a product.

In the rest of the terms at least one $w_j$ is replaced by $B_j$ or $C_j$. Then, since
\begin{eqnarray*}
B_j & = & \lambda\left((u^*x)^{k_j - 1}(x^*u)^{l_j - 1} - \delta_{k_j - 1, l_j - 1}\lambda^{k_j - 1}\right) \\
C_j & = & \gamma\lambda\left((u^*x)^{r_j - 1}(x^*u)^{s_j - 1} - \delta_{r_j - 1,s_j - 1}\lambda^{r_j - 1}\right),
\end{eqnarray*}
we see that such a term is once again
\begin{eqnarray*}
b_0w'_1b_1 \cdots w'_nb_n,
\end{eqnarray*}
so of the same form as $W$ in Equation $(\ref{but})$, but now with the total number of symbols $u^*$ and $x$ strictly smaller than the total number of such symbols in $W$. Thus applying the replacement procedure to each of these terms repeatedly, we finally get $\omega_p(W) = \omega_p(\sum W_i)$, where each $W_i$ has the same form as $W$ in Equation (\ref{but}), but for which the substrings $w_j$ are  either as in Equation $(\ref{word2})$ with $k_j$ or $l_j$ equal to zero, or $w_j$ is as Equation $(\ref{reduction})$ (so that no further replacements can be performed). But then each $W_i$ can be rewritten as a nontrivial alternating product in $\Omega$ and $\{u, u^*\}$, so as before $\omega_p(W_i) = 0$. Thus $\omega_p(W) = 0$.
\end{proof}

We finish at last the proof of Proposition $\ref{comp}$. We know that $pMp = W^*(pqp, v, P)$, and thanks to Lemma \ref{freeness}, $W^*(pqp)$ and $W^*(v, P)$ are $\ast-$free in $pMp$ with respect to $\omega_p$. As $pqp$ is with no atoms with respect to $\omega_p$, with the previous notation, we get $W^*(pqp) \cong L(\Z)$. Concerning $W^*(v, P)$, let
\begin{eqnarray*}
e_{ij} & = & v^i(p - P)(v^*)^j, \\
f_{kl} & = & v^k(P - vv^*)(v^*)^l,
\end{eqnarray*}
for $i, j, k, l \in \N$. With straightforward computations, we see that $(e_{ij})_{i, j \in \N}$ and $(f_{kl})_{k, l \in  \N}$ are systems of matrix units, for all $i, j, k, l \in \N$, $e_{ij}f_{kl} = f_{kl}e_{ij} = 0$, and $W^*(e_{ij}, f_{kl}) = W^*(v, P)$. Moreover, $\omega_p(e_{ii}) = (1 - \gamma)\lambda^i$ and $\omega_{p}(f_{kk}) = (\gamma - \lambda)\lambda^k$, with $\gamma = \beta/(1 - \alpha)$. Consequently, with notation of Proposition \ref{comp}, we finally get $(W^*(v, P), \omega_p) \cong (\C^2, \tau_{\delta}) \otimes (B(\ell^2(\N)), \psi_{\lambda})$. The proof is complete.
\end{proof}
\begin{nota}
\emph{For a von Neumann $(A, \phi_A)$ endowed with a state $\phi_A$, we will denote by $A^\circ$ the kernel of $\phi_A$ on $A$.}
\end{nota}
The next proposition is in some sense a generalization of Theorem $1.2$ of \cite{dykema93}. For the convenience of the reader, we will write a complete proof.
\begin{prop}\label{tensor}
Let $(A, \phi_A)$, $(B, \phi_B)$ and $(C, \phi_C)$ be three von Neumann algebras endowed with faithful, normal states such that $A$ is a factor of type ${\rm I}$.  Let
$${\begin{array}{ccl}
(\mathcal{M}, \psi) & = & \left((C, \phi_C)\otimes(A, \phi_A)\right) \ast (B, \phi_B) \\
\cup & & \\
(\mathcal{N}, \psi) & = & (A, \phi_A) \ast (B, \phi_B) \\
\end{array}}$$
and let $e$ be a minimal projection of $A$. Then in $e\mathcal{M}e$, we have that $e\mathcal{N}e$ and $C \otimes e$ are free with respect to $\psi_e = \frac{1}{\psi(e)}\psi$ and together they generate $e\mathcal{M}e$, so that
$$(e\mathcal{M}e, \psi_e) \cong (C, \phi_C) \ast (e\mathcal{N}e, {\psi_e}).$$
\end{prop}
\begin{proof}
We follow step by step the proof of Theorem $1.2$ of \cite{dykema93}. For notational convenience, we identify $C$ with $C \otimes 1 \subset \mathcal{M}$. To see that $e\mathcal{N}e$ and $eC$ generate $e\mathcal{M}e$, note that $\mathcal{N}$ and $eC$ generate $\mathcal{M}$; so span~$\Lambda(\mathcal{N}, eC)$ is dense in $\mathcal{M}$ and $e\Lambda(\mathcal{N}, eC)e = \Lambda(e\mathcal{N}e, eC)$.

We shall show that $\psi_e$ is zero on a nontrivial alternating product in $(e\mathcal{N}e)^\circ$ and $eC^\circ$. Let $a = e - \psi(e)1$. Then $A^\circ = \C a + S$ where 
$$S = \{s \in A \:|\: \psi(s) = 0, \, ese =~0\}.$$
Let $x \in (e\mathcal{N}e)^\circ$. Then by Kaplansky Density Theorem, $x$ is a s.o.-limit of a bounded sequence $(R_k)_{k \in \N}$ in span~$\Lambda(\{a\} \cup S, B^\circ)$. For $Q \in$ span~$\left(\Lambda(\{a\} \cup S, B^\circ)\backslash\{1\}\right)$, we see that $\psi$ on $eQe$ is equal to a fixed constant times the coefficient of $a$ in $Q$. So since $\psi(R_k) \to 0$ and $\psi(eR_ke) \to 0$, we may assume that the coefficients in each $R_k$ of $1$ and $a$ are zero. Since $R_k - eR_ke \to 0$ for the s.o. topology, we may also assume that the coefficient of each element of $S$ in $R_k$ is zero, i.e., that each $R_k \in$ $\mbox{span} \left(\Lambda(\{a\} \cup S, B^\circ)\backslash(\{1, a\} \cup S)\right)$. To prove the proposition, it suffices to show that $\psi$ is zero on a nontrivial alternating product in $\Lambda(\{a\} \cup S, B^\circ)\backslash(\{1, a\} \cup S)$ and $eC^\circ$. But regrouping and multiplying some neighboring elements gives (a constant times) a nontrivial alternating product in $\{a\} \cup S \cup (eC^\circ) \cup (SC^\circ)$ and $B^\circ$. Thus by freeness, $\psi$ is zero on such a product.
\end{proof}

\begin{proof}[Proof of Theorem $\ref{premierresultat}$]
Apply Proposition $\ref{tensor}$ with $(A, \phi_A) = (B(\ell^2(\N)), \psi_{\lambda})$, $(B, \phi_B) = (L(\Z), \tau)$, $(C, \phi_C) = (\C^2, \tau_{\delta})$. Let $e = e_{00} \in B(\ell^2(\N))$, and denote 
$$\begin{array}{ccl}
(\mathcal{M}, \psi) & = & \left((\C^2, \tau_{\delta}) \otimes (B(\ell^2(\N)), \psi_{\lambda})\right) \ast(L(\Z), \tau) \\
 \cup & & \\
(\mathcal{N}, \psi) & = & (B(\ell^2(\N)), \psi_{\lambda}) \ast (L(\Z), \tau).
\end{array}$$
We get
\begin{equation*}
(e\mathcal{M}e, \psi_e) \cong (\C^2, \tau_{\delta}) \ast (e\mathcal{N}e, \psi_e).
\end{equation*}
But with notation of Section $\ref{bac}$, $(\mathcal{N}, \psi) \cong (T_{\lambda}, \varphi_{\lambda})$ is the free Araki-Woods factor of type ${\rm III_{\lambda}}$. Since $e = e_{00}$, applying Proposition $\ref{compression}$, we get $(e\mathcal{N}e, \psi_e) \cong (T_{\lambda}, \varphi_{\lambda})$. We use now the {\textquotedblleft free absorption\textquotedblright} properties of $(T_{\lambda}, \varphi_{\lambda})$. Denote by $L(\F(s))$ the interpolated free factor with $s$ generators. We know that $(T_{\lambda}, \varphi_{\lambda}) \ast (L(\F_{\infty}), \tau) \cong (T_{\lambda}, \varphi_{\lambda})$ (Corollary $5.5$ in \cite{shlya97}) and $L(\Z) \ast (\C^2, \tau_{\delta}) \cong L(\F(1 + 2\delta(1 - \delta)))$ (Lemma $1.6$ in \cite{dykema93}). Consequently,
\begin{equation*}
(e\mathcal{M}e, \psi_e) \cong (\C^2, \tau_{\delta}) \ast (T_{\lambda}, \varphi_{\lambda}) \cong (T_{\lambda}, \varphi_{\lambda}).
\end{equation*}
But, in a canonical way
\begin{equation*}
(\mathcal{M}, \psi) \cong (e\mathcal{M}e, \psi_e) \otimes (B(\ell^2(\N)), \psi_{\lambda}).
\end{equation*}
Since $(T_{\lambda}, \varphi_{\lambda}) \cong (T_{\lambda}, \varphi_{\lambda}) \otimes (B(\ell^2(\N)), \psi_{\lambda})$, we get 
\begin{equation*}
(\mathcal{M}, \psi) \cong (T_{\lambda}, \psi_{\lambda}).
\end{equation*}
We remind that we have proved (Proposition $\ref{comp}$) that 
\begin{equation*}
\left(pMp, \omega_p \right)  \cong  L(\Z) \ast \left((\C^2, \tau_{\delta}) \otimes (B(\ell^2(\N)), \psi_{\lambda})\right) = (\mathcal{M}, \psi),
\end{equation*}
where $(M, \omega) = (M_2(\C), \omega_{\lambda}) \ast (\C^2, \tau_{\beta})$ and $p = p_{11} \in M_2(\C)$. Thus,
\begin{equation*}
\left(pMp, \omega_p\right)  \cong  (T_{\lambda}, \varphi_{\lambda}).
\end{equation*}
But once again 
\begin{equation*}
(M, \omega) \cong (pMp, \omega_p) \otimes (M_2(\C), \omega_{\lambda}).
\end{equation*}
Since $(T_{\lambda}, \varphi_{\lambda}) \cong (T_{\lambda}, \varphi_{\lambda}) \otimes (M_2(\C), \omega_{\lambda})$, we finally get 
\begin{equation*}
(M_2(\C), \omega_{\lambda}) \ast (\C^2, \tau_{\beta}) \cong (T_{\lambda}, \varphi_{\lambda}).
\end{equation*}
\end{proof}

\section{Proof of the Main Theorem}\label{mainthe}

The aim of this section is to prove Theorem $\ref{results}$. We will be using the {\textquotedblleft machinery\textquotedblright} of amalgamated free products of von Neumann algebras. We introduce some notations and recall some basic facts about free products with amalgamation (see \cite{{popa1993}, {ueda}, {voiculescu85}}).

Let $(B, \varphi_B)$, $(A_i, \varphi_i)$, $i = 1, 2$, be three von Neumann algebras endowed with faithful normal states. Assume that there exist modular embeddings $\rho_i : (B, \varphi_B) \hookrightarrow (A_i, \varphi_i)$. Denote by $E_i : A_i \to B$ the unique state-preserving conditional expectation associated with the embedding $\rho_i$. We shall denote by
\begin{equation*}
(M, E) := (A_1, E_1) \mathop{\ast}_B (A_2, E_2)
\end{equation*}
the free product with amalgamation over $B$ of $A_1$ and $A_2$ w.r.t. the conditional expectations $E_1$ and $E_2$.

Let $(B, \varphi_B)$ and $(C, \varphi_C)$ be two von Neumann algebras together with a faithful normal state. Let $(A, \varphi_A) = (B, \varphi_B) \ast (C, \varphi_C)$ be their free product. We have canonical modular embeddings $\rho_B : (B, \varphi_B) \hookrightarrow (A, \varphi_A)$, $\rho_C : (C, \varphi_C) \hookrightarrow (A, \varphi_A)$ (see \cite{dykema94b} for further details). We shall regard $B, C \subset A$. Define as before $F : A \to B$ to be the (unique) state-preserving conditional expectation. Let $B^\circ = B \cap \ker(\varphi_B)$, $C^\circ = C \cap \ker(\varphi_C)$ and denote as usual $\Omega = \Lambda(B^\circ, C^\circ)$ the set of alternating products in $B^\circ$ and $C^{\circ}$ including the trivial one. From \cite{dykema96}, we know that 
\begin{eqnarray*}
\forall b \in B, F(b) & = & b \\
\forall z \in \Omega \backslash (B^\circ \cup \{1\}), F(z) & = & 0.
\end{eqnarray*}
The following proposition is well-known from specialists but we will give a proof for the sake of completeness.
\begin{prop}\label{amal}
We use the same notations as before. Moreover, let $(M, \varphi_M)$ be a von Neumann algebra such that $B \subset M$ together with $E : (M, \varphi_M) \to (B, \varphi_B)$ a state-preserving conditional expectation.
Denote by $(\mathcal{M}, G) = (M, E) \mathop{\ast}_{B} \left(B \ast C, F\right)$ and denote by $\psi = \varphi_B \circ G$ the canonical state on $\mathcal{M}$. Then,
\begin{equation*}
(\mathcal{M}, \psi) \cong (M, \varphi_M) \mathop{\ast} (C, \varphi_C).
\end{equation*}
\end{prop}

\begin{proof}
We see immediatly that $(M, \varphi_M)$ and $(C, \varphi_C)$ embed in $(\mathcal{M}, \psi)$ in a state-preserving way and together they generate $\mathcal{M}$. It remains to prove that $M$ and $C$ are free together w.r.t. the state $\psi$. For notational convenience, we may assume $M, C \subset \mathcal{M}$. Denote $M^{\circ} = M \cap \ker(\varphi_M)$, $C^{\circ} = C \cap \ker(\varphi_C)$. Let $W$ be a nontrivial alternating product in $M^{\circ}$ and $C^{\circ}$, so that $W$ can be written
$$W = x_0w_1x_1 \cdots w_nx_n,$$
where $x_j \in C^{\circ}$, $w_j \in M^{\circ}$ for all $j$, except possibly $x_0$ and/or $x_n$ are equal to $1$. Denote $\Omega = \Lambda(B^\circ, C^{\circ})$. If $W \in M^{\circ}$, there is nothing to do. If not, for each $j$, replace $w_j$ by
\begin{equation*}
w_j = w'_j + b_j,
\end{equation*}
where $w'_j \in M \cap \ker E$ and $b_j \in B^\circ$. Applying the replacement procedure and multiplying some neighboring elements, we get $\psi(W) = \psi(\sum W_i)$ where each $W_i$ is a nontrivial alternating product in $M \cap \ker E$ and $\Omega \backslash (B^\circ \cup \{1\})$. But, we saw that $\Omega \backslash (B^\circ \cup \{1\}) \subset \ker F$. Thus, by freeness with amalgamation over $B$, we get $G(W_i) = 0$. But $\psi(W_i) = (\psi \circ G)(W_i) = 0$, consequently $\psi(W) = 0$.
\end{proof}

\begin{lem}\label{absorb}
Let $(N, \varphi)$ be a von Neumann algebra endowed with a faithful normal state such that the centralizer $N^\varphi$ is a factor. Let $(M_n(\C), \omega)$ be a matrix algebra endowed with a faithful normal state. Let $\rho_i : (M_n(\C), \omega) \hookrightarrow (N, \varphi)$, $i = 1, 2$, be two modular embeddings. Then, there exists a unitary $u \in  \mathcal{U}(N^\varphi)$ such that $\Ad(u) \circ \rho_1 = \rho_2$. 
\end{lem}
\begin{proof}
Denote by $(e_{kl})_{0 \leq k, l \leq n - 1}$ the matrix unit in $M_n(\C)$.
Let $i \in \{1, 2\}$. Denote $p_i = \rho_i(e_{00})$. Since $\rho_i$ is modular, we have $p_{1, 2} \in N^\varphi$ and $\varphi(p_1) = \varphi(p_2) = \omega(e_{00})$. Since $N^\varphi$ is a factor, there exists a partial isometry $v \in N^\varphi$ such that $p_1 = v^*v$ and $p_2 = vv^*$. Denote $u = \sum_{i = 0}^{n - 1} \rho_2(e_{i0}) v \rho_1(e_{0i})$. An easy computation shows that $u$ is a unitary and $u \in N^\varphi$, since $\rho_{1, 2}$ are modular. Moreover, for any $0 \leq k, l \leq n - 1$, 
\begin{equation*}
u \rho_1(e_{kl}) u^* = \rho_2(e_{kl}). 
\end{equation*}
\end{proof}

\begin{theo}\label{lemtheo}
Let $(A_1, \phi_1)$ and $(A_2, \phi_2)$ be any von Neumann algebras endowed with faithful, normal states. Assume that for some $\lambda, \beta \in ]0, 1[$, there exist modular embeddings
\begin{eqnarray*}
(M_2(\C), \omega_\lambda) & \hookrightarrow & (A_1, \phi_1) \\
(\C^2, \tau_\beta) & \hookrightarrow & (A_2, \phi_2),
\end{eqnarray*}
such that $\lambda/(\lambda + 1) \leq \min\{\beta, 1 - \beta\}$. Then
\begin{equation*}
(A_1,Ê\phi_1) \ast (A_2, \phi_2) \cong (A_1, \phi_1) \ast (T_\lambda, \varphi_\lambda) \ast (A_2, \phi_2).
\end{equation*}
\end{theo}

\begin{proof}
We shall simply denote by $M_2$ the matrix algebra $M_2(\C)$. Denote by $(A, \phi) = (A_1, \phi_1) \ast (A_2, \phi_2)$ and denote by
\begin{eqnarray*}
E_1 : (A_1, \phi_1) & \to & (M_2, \omega_\lambda) \\
E_2 : (A_2, \phi_2) & \to & (\C^2, \tau_\beta) \\
\widetilde{E} : (A, \phi) & \to & (\C^2, \tau_\beta)
\end{eqnarray*}
the canonical state-preserving conditional expectations. Since $(M_2, \omega_{\lambda}) \ast (\C^2, \tau_\beta) \cong (T_\lambda, \varphi_\lambda)$ (Theorem $\ref{premierresultat}$), denote by $F : (T_\lambda, \varphi_\lambda) \to (M_2, \omega_\lambda)$ the associated state-preserving conditional expectation. Let $(N, \psi) = (A_1,\phi_1) \ast (\C^2, \tau_\beta)$.
Applying Proposition $\ref{amal}$, we get
\begin{equation*}
(N, E) \cong \displaystyle{(A_1,E_1) \mathop{\ast}_{M_2}(M_2 \ast \C^2, F)}.
\end{equation*}
Since $(T_\lambda)^{\varphi_\lambda} \cong L(\mathbf{F}_\infty)$ is a factor, applying Lemma $\ref{absorb}$ for $n = 2$, we obtain that the modular embedding of $(M_2, \omega_\lambda)$ into $(T_\lambda, \varphi_\lambda)$ is unique up to a conjugation by a unitary in $T_{\lambda}^{\varphi_\lambda}$, and we have
\begin{eqnarray*}
(N, E) & \cong & \displaystyle{(A_1,E_1) \mathop{\ast}_{M_2} (M_2 \ast T_{\lambda}, F)} \\
(N, \psi) & \cong & (A_1, \phi_1) \ast (T_\lambda, \varphi_\lambda) \; (\mbox{by Proposition } \ref{amal}) \\
& \cong & \left((A_1, \phi_1) \ast (T_\lambda, \varphi_\lambda)\right) \ast (T_\lambda, \varphi_\lambda).
\end{eqnarray*}
From Theorem $11$ of \cite{barnett95}, we get that the centralizer algebra $N^\psi$ is a factor. If $\rho_i : (\C^2, \tau_\beta) \hookrightarrow (N, \psi)$ are two modular embeddings, denote $p_i = \rho_i(p) \in N^\psi$ such that $\psi(p_i) = \beta$. Since $p_1$ and $p_2$ are equivalent in $N^\psi$, $\rho_1$ and $\rho_2$ are unitarily conjugate. Consequently, using the isomorphism 
\begin{equation}\label{neweq}
(N, \psi)  \cong  (A_1, \phi_1) \ast (T_\lambda, \varphi_\lambda) \ast (\C^2, \tau_\beta),
\end{equation}
we shall denote by $G : (N, \psi) \to (\C^2, \tau_\beta)$ the associated state-preserving conditional expectation. We finally get
\begin{eqnarray*}
(A, \widetilde{E}) & \cong & \left(A_1 \ast \C^2, G\right) \mathop{\ast}_{\C^2} (A_2, E_2) \; (\mbox{by Proposition } \ref{amal})\\
& \cong & \left(A_1 \ast T_\lambda \ast \C^2, G\right) \mathop{\ast}_{\C^2} (A_2, E_2) \; (\mbox{by } (\ref{neweq}))\\
(A, \phi) & \cong & (A_1, \phi_1) \ast (T_\lambda, \varphi_\lambda) \ast (A_2, \phi_2). \; (\mbox{by Proposition } \ref{amal})
\end{eqnarray*}
The proof is complete.
\end{proof}

Theorem $\ref{results}$ is a straightforward corollary of Theorem $\ref{lemtheo}$. We end this section by giving some examples of von Neumann algebras which satisfy assumptions of Theorem $\ref{results}$. We introduce  the class $\mathcal{S}$ of all von Neumann algebras $(M, \phi)$ with separable predual and endowed with a faithful, normal, almost periodic state $\phi$ such that
\begin{equation*}
(M, \phi) \ast (L(\Z), \tau_{\Z}) \cong (T_{\Sd(\phi)}, \varphi_{\Sd(\phi)}).
\end{equation*}
Note that if $(M_1, \phi_1)$ and $(M_2, \phi_2)$ are in $\mathcal{S}$, then $(M_1, \phi_1) \ast (M_2, \phi_2)$ is also in $\mathcal{S}$.
\begin{exam}\label{examabs}
\emph{We give several examples of von Neumann algebras in the class $\mathcal{S}$. This list is not exhaustive and there is nothing really new here: these examples are mere consequences of results in \cite{{dykema96}, {dykema94}, {dykema93}, {radulescu1994}, {shlya97}, {voiculescu92}}, and of Proposition $\ref{tensor}$.}
\begin{enumerate}

\item [$\bullet$] \emph{Type ${\rm I}$: All factors of type ${\rm I}$ endowed with a faithful, normal nontracial state~$\phi$.}

\item [$\bullet$] \emph{Type ${\rm III}$: All the almost periodic free Araki-Woods factors $(T_{\Gamma}, \varphi_{\Gamma})$ endowed with their free quasi-free state.}

\item [$\bullet$] \emph{Tensor products: All the tensor products $(N, \omega) \otimes (\mbox{Type } {\rm I}, \phi)$, where $(N, \omega)$ is:
\begin{enumerate}
\item [$-$] Any finite-dimensional von Neumann algebra of the form $\displaystyle{\mathop{\C}_{\alpha_1} \oplus \cdots \oplus \mathop{\C}_{\alpha_n}}$ with $\alpha_i >~0$ for all $i$ and $\sum \alpha_i = 1$.
\item [$-$] $(\mathcal{R}, \tau)$ the hyperfinite ${\rm II_1}$ factor.
\item [$-$] Any interpolated free group factor $L(\F(s))$, $s > 1$.
\end{enumerate} }

\item [$\bullet$] \emph{Free products: All the free products of the previous examples.}
\end{enumerate}
\end{exam}

\section{Remark}\label{rema}

We still do not know whether all the free products of finite dimensional matrix algebras $(A_1, \phi_1) \ast (A_2, \phi_2)$ are isomorphic to free Araki-Woods factors. Assume that $A_1 = M_n(\C)$ with
\begin{equation*}
\phi_1 = \Tr\left(
\begin{pmatrix}
\lambda_1 & & \\
& \ddots & \\
& & \lambda_n
\end{pmatrix} \cdot
\right), \: \lambda_1 \leq \cdots \leq \lambda_n.
\end{equation*}
Let $\beta \in ]0, 1[$ such that $\lambda_1 \leq \min\{\beta, 1 - \beta\}$. With our techniques, it is not difficult to see that if one can prove that $(A_1, \phi_1) \ast (\C^2, \tau_{\beta})$ is a free Araki-Woods factor, then all the free products $(A_1, \phi_1) \ast (A_2, \phi_2)$ are also free Araki-Woods factors. That is exactly what we did in Section $3$ for $n = 2$. But one of the crucial ingredients in the proof was the precise picture of Voiculescu in Theorem $\ref{deuxprojs}$. This precise description no longer exists for $n \geq 3$ (see \cite{dykema93} for further details).

\begin{ack}
\emph{The author would like to thank Stefaan Vaes for helpful suggestions regarding this manuscript.}
\end{ack}

\bibliographystyle{plain}

\begin{thebibliography}{AA}


\bibitem{barnett95} {\sc L. Barnett}, {\it Free product von Neumann algebras of type ${\rm III}$}. Proc. Amer. Math. Soc. {\bf123} (1995), 543--553.

\bibitem{connes74} {\sc A. Connes}, {\it Almost periodic states and factors of type ${\rm III_1}$}. J. Funct. Anal. {\bf16} (1974), 415--445.


\bibitem{connes73} {\sc A. Connes},
{\it Une classification des facteurs de type {\rm III}.} Ann. Sci. {\'E}cole Norm. Sup. {\bf 6} (1973), 133--252.


\bibitem{dykema96} {\sc K. Dykema},
{ \it Free products of finite-dimensional and other von Neumann algebras
with respect to non-tracial states.} Free probability theory (Waterloo, ON, 1995)
Fields Inst. Commun. {\bf 12}
Amer. Math. Soc., Providence, RI, 1997, pp. 41--88.

\bibitem{dykema94} {\sc K. Dykema}, {\it Interpolated free group factors.} Pacific J. Math. {\bf 163} (1994), 123--135.

\bibitem{dykema94b} {\sc K. Dykema}, {\it Factoriality and Connes' invariant $T(M)$ for free products of von Neumann algebras.} J. reine angew. Math. {\bf 450} (1994), 159--180.

\bibitem{dykema93} {\sc K. Dykema}, {\it Free products of hyperfinite von Neumann algebras and free dimension.} Duke Math. J. {\bf 69} (1993), 97--119.

%


\bibitem{popa1993} {\sc S. Popa}, {\it Markov traces on universal Jones algebras and subfactors of finite index.} Invent. Math. {\bf 111}, (1993), 375--405.

\bibitem{radulescu1994} {\sc F. R\u{a}dulescu}, { \it Random matrices, amalgamated free products and subfactors of the von Neumann algebra of a free group, of noninteger index.} Invent. Math. {\bf 115}, (1994), 347--389.

\bibitem{radulescu1993} {\sc F. R\u{a}dulescu}, {\it A type ${\rm III_{\lambda}}$ factor with core isomorphic to the von Neumann algebra of a free group, tensor $B(H)$.} Recent Advances in Operator Algebras (Orl\'eans, 1992). Ast\'erisque {\bf 232} (1995), 203--209.





\bibitem{shlya2004} {\sc D. Shlyakhtenko}, {\it On the Classification of Full
  Factors of Type {\rm III}.} Trans. Amer. Math. Soc. {\bf 356} (2004), 4143--4159.

\bibitem{shlya99} {\sc D. Shlyakhtenko}, {\it $A$-valued semicircular systems.}
    J. Funct. Anal. {\bf 166} (1999), 1--47.

\bibitem{shlya98} {\sc D. Shlyakhtenko}, {\it Some applications of freeness with amalgamation.} J. Reine Angew. Math. {\bf 500} (1998), 191--212.

\bibitem{shlya97} {\sc D. Shlyakhtenko}, {\it Free quasi-free states.} Pacific J. Math. {\bf 177} (1997), 329--368.

\bibitem{takesakiII} {\sc M. Takesaki}, { \it Theory of Operator Algebras ${\rm II}$.} {\it EMS} {\bf 125}. Springer-Verlag, Berlin, Heidelberg, New-York, 2000.

\bibitem{ueda} {\sc Y. Ueda}, {\it Amalgamated free products over Cartan subalgebra.} Pacific J. Math. {\bf 191} (1999), 359--392.



\bibitem{vaes2004} {\sc S. Vaes}, {\it \'Etats quasi-libres libres et facteurs de
  type {\rm III} (d'apr{\`e}s D. Shlyakhtenko).} S{\'e}minaire Bourbaki, expos\'e 937, {Ast\'erisque {\bf
  299}} (2005), 329--350.


\bibitem{voiculescu92} {\sc D.-V. Voiculescu, K.J. Dykema \& A. Nica}, {\it Free
  random variables.} CRM Monograph Series {\bf 1}.
American Mathematical Society, Providence, RI, 1992. 

\bibitem{voiculescu90} {\sc D.-V. Voiculescu}, {\it Circular and semicircular systems and free product factors.} Operator algebras, Unitary representations, Enveloping Algebras, and Invariant Theory, Progress in Mathematics {\bf 92}.
Birkh\"{a}user, Boston, (1990), 45--60. 

\bibitem{voiculescu87} {\sc D.-V. Voiculescu}, {\it Multiplication of certain noncommuting random variables.} J. Operator Theory {\bf 18} (1987), 223--235.

\bibitem{voiculescu85} {\sc D.-V. Voiculescu}, {\it Symmetries of some reduced free product $C^*$-algebras.} Operator algebras and Their Connections with Topology and Ergodic Theory, Lecture Notes in Mathematics {\bf 1132}.
Spriger-Verlag, (1985), 556--588. 

\end{thebibliography}

\end{document}